\documentclass[a4paper,12pt]{article}
\usepackage{amsmath,amsfonts,amsthm} 
\usepackage{graphicx}

\newcommand{\B}[1]{\boldsymbol{#1}}

\newcommand{\C}{{\mathbb C}}

\def\proof{\medskip\noindent {\it Proof.}\quad }
\def\Mat{{\rm Mat}}
\newcommand{\G}{\Gamma}

\newtheorem{theorem}{Theorem}[section]
\newtheorem{lemma}[theorem]{Lemma}
\newtheorem{corollary}[theorem]{Corollary}

\theoremstyle{definition}
\newtheorem{note}[theorem]{Note}
\newcommand{\de}{\em}

\begin{document}
\pagestyle{plain}

\title{A characterization of $Q$-polynomial \\
            distance-regular graphs}
\author{Aleksandar Juri\v{s}i\'{c}, Paul Terwilliger, Arjana \v{Z}itnik}
\date{August 3, 2009}
\maketitle

\begin{abstract} \noindent
%
We obtain the following characterization of
$Q$-polynomial distance-regular graphs.
Let $\G$ denote a distance-regular graph with diameter $d\ge 3$.
Let $E$ denote a minimal idempotent of $\G$
which is not the trivial idempotent $E_0$.
Let $\{\theta_i^*\}_{i=0}^d$
denote the dual eigenvalue sequence for $E$.
We show that $E$ is $Q$-polynomial if and only if 
(i)   the entry-wise product $E \circ E$ is a linear combination 
      of $E_0$, $E$, and at most one other minimal idempotent of $\G$;
(ii)  there exists a complex scalar $\beta$ such that
      $\theta^*_{i-1}-\beta \theta^*_i + \theta^*_{i+1}$
		is independent of $i$ for $1 \le i \le d-1$;
(iii) $\theta^*_i \ne \theta^*_0$ for $1 \le i \le d$.
\end{abstract}

\normalbaselineskip=18pt
\baselineskip=\normalbaselineskip

\section{Introduction}

In this paper we give a new characterization of the $Q$-polynomial 
property for distance-regular graphs. 
In order to motivate and describe our result, we first recall some notions. 
Let $\G$ denote a distance-regular graph with diameter $d \ge 3$
and vertex set $X$ (see Section 2 for formal definitions).  
Recall that there exists a minimal idempotent $E_0$ of $\G$ 
such that $E_0=|X|^{-1}J$, where $J$ is the all 1's matrix. 
We call $E_0$ {\it trivial}.
Let $\{ E_i\}_{i=1}^d$ denote an ordering of the nontrivial minimal 
idempotents of $\G$.
It is known that for $0\le i,j\le d$ the entry-wise product 
$E_i \circ E_j$ is a linear combination of the minimal idempotents 
of $\G$, so that
$$ %
E_i \circ E_j \ =\ {1\over |X|}\, \sum_{h=0}^{d}\, q_{ij}^h\, E_h. 
$$ %
The coefficients $q^h_{ij}$ are called the {\it Krein parameters} of $\G$. 
They are real and nonnegative; see for example \cite[p.~48--49]{BCN}.
%
Now consider when is a Krein parameter zero. 
Note that $J \circ E_j = E_j$ for $0\le j\le d$, 
so $q^h_{0j} = \delta_{hj}$   for $0\le h,j\le d$.
The ordering $\{ E_i\}_{i=1}^d$ is called {\it $Q$-polynomial} whenever
for all $0\le i,j\le d$ the Krein parameter $q^1_{ij}$ is zero 
if $|i-j|>1$ and nonzero if $|i-j|=1$.
Let $E$ denote a nontrivial minimial idempotent of $\G$.
We say that $E$  is {\it $Q$-polynomial} 
whenever there exists a $Q$-polynomial ordering $\{ E_i\}_{i=1}^d$ 
of the nontrivial minimal idempotents of $\G$ such that $E=E_1$. 
We now explain the $Q$-polynomial property 
in terms of representation diagrams. 
Let $E$ denote a nontrivial minimal idempotent of $\G$,
and for notational convenience write $E=E_1$.
The {\it representation diagram} $\Delta_E$
is the undirected graph with vertex set $\{0,\ldots, d\}$ such that 
vertices $i,j$ are adjacent whenever $i \ne j$ and $q^1_{ij} \ne 0$.
By our earlier comments, $q^1_{0j}=\delta_{1j}$ for $0\le j\le d$.
Therefore, in $\Delta_E$ the vertex $0$ is adjacent to the vertex $1$ 
and no other vertex, see Figure 1(a). 
Observe that $E$ is $Q$-polynomial if and only if $\Delta_E$ is a path, 
and in this case the natural ordering $0,1,\ldots$ of the vertices 
in $\Delta_E$ agrees with the
$Q$-polynomial ordering associated with~$E$. See Figure 1(c).

\bigskip
\includegraphics[width=150mm]{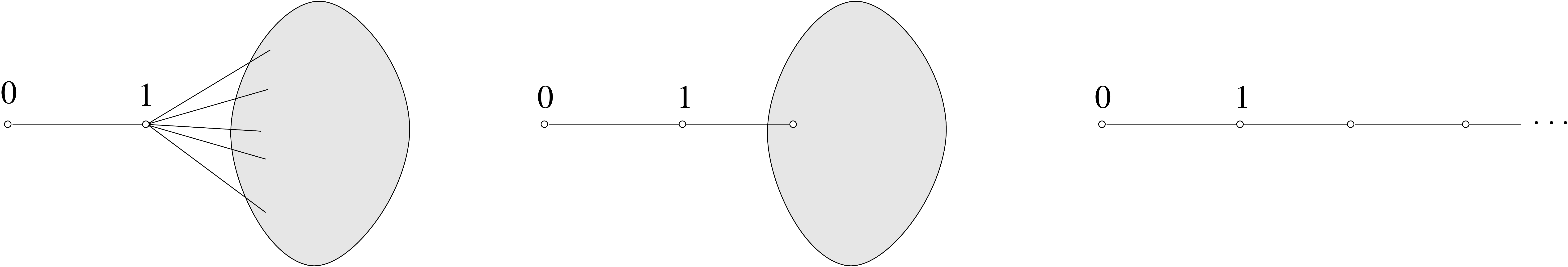}\\
\hspace*{2cm} (a) \hspace*{4cm} (b) \hspace*{5.2cm} (c)
\bigskip

\baselineskip 9pt \par{\leftskip 11mm \rightskip 11mm \noindent
{\footnotesize Figure 1:
For $E=E_1$ the representation diagram $\Delta_E$ around 
vertices 0 and 1:
(a) in general; 
(b) when $E=E_1$ is a tail;
(c) when $E$ is $Q$-polynomial.
}\par} \baselineskip=\normalbaselineskip

\medskip\noindent

Let $E$ denote a nontrivial minimal idempotent of $\G$.
In this paper we give a condition that is necessary and
sufficient for $E$ to be $Q$-polynomial.
We now describe the condition, which has three parts.

The first part has to do with the representation diagram $\Delta_E$.
According to Lang \cite{lang1}, $E$ is a {\it tail} whenever
$E \circ E$ is a linear combination of $E_0$, $E$, 
and at most one other minimal idempotent of $\G$. 
In terms of the diagram $\Delta_E$, and writing $E=E_1$ for notational
convenience, $E$ is a tail if and only if vertex $1$ is adjacent
to at most one vertex besides vertex $0$, see Figure 1(b).
Note that if $E$ is $Q$-polynomial, then $E$ is a tail.

The next part of our condition involves the dual eigenvalue sequence 
$\{ \theta_i^*\}_{i=0}^d$ for E.
This sequence satisfies $E=|X|^{-1}\sum_{i=0}^d \theta_i^* A_i$,
where $\{ A_i\}_{i=0}^d$ are the distance matrices of $\G$. 
Following \cite{lang1}, 
we say that $E$ is {\it three-term recurrent} (in short TTR)
whenever there exists a complex scalar $\beta$ such that
$\theta^*_{i-1}- \beta \theta^*_i + \theta^*_{i+1}$ is independent of
$i$ for $1 \le i \le d-1$. 
If $E$ is $Q$-polynomial, then $E$ is TTR
by \cite[Theorem 8.1.2]{BCN}, cf. \cite{L}.

The third part of our condition involves both the
diagram $\Delta_E$ and the dual eigenvalues 
$\{ \theta_i^*\}_{i=0}^d$.
By \cite[Proposition 2.11.1]{BCN}, $\Delta_E$ is connected 
if and only if $\theta_i^* \ne \theta_0^*$ for $1 \leq i \leq d$. 
These equivalent statements hold if $E$ is $Q$-polynomial,
since in this case $\Delta_E$ is a path.

\bigskip
We now state our main result.

\begin{theorem}                                           \label{main}
Let $\G$ denote a distance-regular graph with diameter $d \ge 3$.
Let $E$ denote a nontrivial minimal idempotent for $\G$ 
and let $\{ \theta_i^*\}_{i=0}^d$
denote the corresponding dual eigenvalue sequence.
Then $E$ is $Q$-polynomial if and only if 
\begin{description} \itemsep 35pt \parskip -30pt
\item[{\rm (i)}]   \hspace*{9pt} $E$ is a tail,
\item[{\rm (ii)}]  \hspace*{6pt}  $E$ is TTR,
\item[{\rm (iii)}] \hspace*{3pt}  $\theta^*_i \ne \theta^*_0$ for $1\le i\le d$.
\end{description}
\end{theorem}

\bigskip
In Section \ref{cx} 
we discuss the minimality of assumptions (i)--(iii) of Theorem \ref{main}.
We show that in general, no proper subset of (i)--(iii) is sufficient 
to imply that $E$ is $Q$-polynomial.
Theorem \ref{main} gives a characterization of the $Q$-polynomial
distance-regular graphs. A similar characterization,
where assumption (i) is replaced by some equations
involving the dual eigenvalues and intersection numbers,
is given by Pascasio \cite{pascasio3}.

\section{Preliminaries}

In this section we review some definitions and basic concepts. 
See Brouwer, Cohen and Neumaier \cite{BCN} and Terwilliger \cite{Talg}
for more background information.
Let $\C$ denote the complex number field and $X$ a nonempty finite set.
Let $\Mat_X\C$ denote the $\C$-algebra consisting of all matrices whose 
rows and columns are indexed by $X$ and whose entries are in $\C$. 
Let $V = \C^X$ denote the vector space over $\C$ consisting of column 
vectors whose coordinates are indexed by $X$ and whose entries are in $\C$.
Observe that $\Mat_X\C$ acts on $V$ by left multiplication.
For all $y \in X$, let $\hat{y}$ denote the element of 
$V$ with $1$ in the $y$-th coordinate and $0$ in all other coordinates. 

From now on $\G$ denotes a finite, undirected, connected
graph, without loops or multiple edges, with vertex set $X$, 
the shortest path-length distance function $\partial$ and
diameter $d:=\max \{\partial(x,y)\, |\, x,y \in X\}$.
For a vertex $x \in X$ and integer  $i \ge 0$
define $$\G_i(x)= \{ y \in X\, |\, \partial (x,y)=i \}.$$
For notational convenience abbreviate $\G(x)=\G_1(x)$.
For an integer $k \ge 0$, the graph $\G$ is said to be 
{\it regular with valency k} whenever $|\G(x)|= k$ for all $x\in X$.
The graph $\G$ is said to be {\de distance-regular} whenever for all integers 
$h,i,j$ $(0 \leq h,i,j\leq d)$ and  vertices $x,y \in X$ with $\partial(x,y)=h$,
the number $p_{ij}^h=|\G_i(x)\cap \G_j(y)|$ is independent of $x,y$.
The constants $p_{ij}^h$ 
are called the {\de intersection numbers} of $\G$.
From now on  assume $\G$ is distance-regular with diameter $d \ge 3$.
Note that $\G$ is regular with valency $k=p^0_{11}$.

We recall the Bose-Mesner algebra of $\G$.
For $0\le i\le d$ let $A_i$ denote the matrix in $\Mat_X\C$ with 
$(x,y)$-entry equal to $1$ if $\partial(x,y)=i$ and $0$ otherwise.
We call $A_i$ the {\de $i$-th distance matrix} of $\G$. 
Note that $A_i$ is real and symmetric. 
We observe that $A_0=I$, where $I$ is the identity matrix,
and abbreviate $A=A_1$. 
We observe that  $\sum_{i=0}^d A_i=J$ and
$A_iA_j=\sum_{h=0}^d p_{ij}^h A_h $ for $0 \le i,j \le d$.
Let $M$ denote the subalgebra of $\Mat_X\C$ generated by $A$.
By \cite[p.~127]{BCN} the matrices $\{ A_i\}_{i=0}^d$ form a basis for $M$. 
We call $M$ the {\de Bose-Mesner algebra} of $\G$. 
%
By \cite[p.~45]{BCN}, $M$ has a basis $\{E_i\}_{i=0}^d$
such that (i) $E_0=|X|^{-1}J$; (ii) $I=\sum_{i=0}^d E_i$;
(iii) $E_iE_j =\delta_{ij}E_i$ for $0 \le i,j\le d$.
%
%
By \cite[p.~59,\, 64]{BI} 
the matrices $\{E_i\}_{i=0}^d$ are real 
and symmetric.
We call $\{ E_i\}_{i=0}^d$ the {\de minimal idempotents} of $\G$.
We call $E_0$ {\it trivial}.
Since $\{ E_i\}_{i=0}^d$ form a basis for $M$, there
exist complex scalars $\{ \theta_i\}_{i=0}^d$ such that 
\begin{equation}                   \label{prodAEi0}
A\ =\ \sum_{i=0}^d \theta_i E_i.
\end{equation}
By (\ref{prodAEi0}) and since $E_i E_j = \delta_{ij}E_i$ we have
\begin{equation}                   \label{prodAEi}
AE_i\, =\, E_iA\, =\, \theta_i E_i \ \ \ \   
\mbox{($0\le i\le d$)}.
\end{equation}
We call the scalar $\theta_i$ the {\de eigenvalue} of $\G$
corresponding to $E_i$. Note that the eigenvalues $\{ \theta_i\}_{i=0}^d$ 
are mutually distinct since $A$ generates $M$. 
Moreover $\{\theta_i\}_{i=0}^d$ are real, since $A$ and 
$\{ E_i\}_{i=0}^d$ are real. 
Let $E$ denote a minimal idempotent of $\G$.
Since $\{ A_i\}_{i=0}^d$ form a basis for $M$, there
exist  complex scalars $\{ \theta_i^*\}_{i=0}^d$ such that
\begin{equation}                             \label{dualeigenvaluesdef}
E\ =\ {1\over |X|}\, \sum_{i=0}^d\, \theta_i^*\, A_i.
\end{equation}
We call $\theta_i^*$ the {\de $i$-th dual eigenvalue} of $\G$
corresponding to $E$.  Note that $\{\theta^*_i\}_{i=0}^d$ are real,
since $E$ and $\{ A_i\}_{i=0}^d$ are real.
Let $\circ$ denote the entry-wise product in $\Mat_X\C$.
Observe $A_i \circ A_j =\delta_{ij}A_i$ for $0\le i,j\le d$,
so $M$ is closed under $\circ$.
%
%
Therefore there exist complex scalars $q^h_{ij}$ such that
\begin{equation}\label{kreindef}
E_i \circ E_j = {1\over |X|}\, \sum_{h=0}^d q^h_{ij} E_h
\ \ \ \ \ \ \ (0 \le i,j\le d).
\end{equation}
We call the $q^h_{ij}$ the {\de Krein parameters} of $\G$.
These parameters are real and nonnegative \cite[p.~48--49]{BCN}.
For the moment fix integers $h,i,j$ $(0 \leq h,i,j\leq d)$.
By construction $q^h_{ij}=q^h_{ji}$. 
By \cite[Lemma 2.3.1]{BCN} we have 
$m_h q^h_{ij} = m_i q^i_{jh} = m_j q^j_{hi}$.
Therefore
\begin{eqnarray}                                        \label{kreinsym}
          q^h_{ij}=0  \ \ \ \ \mbox{iff}  \ \ \ \
          q^i_{jh}=0 \ \ \  \ \mbox{iff}  \ \ \ \
          q^j_{hi}=0.  
\end{eqnarray}
%

We recall the dual Bose-Mesner algebra of $\G$ \cite[p.~378]{Talg}. 
For the rest of this section, fix a vertex $x \in X$. 
For $0\le i \le d$ let $E_i^* = E_i^*(x)$ 
denote the diagonal matrix in $\Mat_X\C$ with $(y,y)$-entry
$$
(E_i^*)_{yy}\, =\, \left\{ \begin{array}{ll}
                          1; & \mbox{if} \ \ \partial(x,y)=i\\
                          0; & \mbox{if} \ \ \partial(x,y)\ne i\\
                           \end{array} \right. 
\ \ \ \ \ (y \in X).
$$
%
%
We call $E_i^*$ the {\de $i$-th  dual idempotent of $\G$ with respect to $x$.}
We observe  that $I=\sum_{i=0}^d E_i^*$
and  $E_i^*E_j^*=\delta_{ij}E_i^*$ for $0\le i,j\le d$.
Therefore the matrices $\{ E_i^*\}_{i=0}^d$ form a basis for 
a commutative subalgebra $M^*=M^*(x)$ of  $\Mat_X\C$.
We call $M^*$ the {\de dual Bose-Mesner algebra} of $\G$
with respect to $x$. 

For $0\le i \le d$ let $A_i^*=A_i^*(x)$ denote the diagonal matrix
in $\Mat_X\C$ with $(y,y)$-entry
$$
(A_i^*)_{yy}\, =\, |X|\,(E_i)_{xy} \ \ \ \ \ (y \in X).
$$
We call $A_i^*$ the {\it dual distance matrix} corresponding to $E_i$.
By \cite[p.~379]{Talg} the matrices $\{ A_i^*\}_{i=0}^d$ form a 
basis for $M^*$.
Select an integer $\ell$ $(1\le \ell \le d)$ and
set $E := E_\ell$, $A^* := A_\ell^*$. 
Let  $\{ \theta_i^*\}_{i=0}^d$ denote the dual eigenvalues 
corresponding to $E$. Then using (\ref{dualeigenvaluesdef}) 
we find
\begin{equation}   \label{productdual01}
A^*\ =\ \sum_{i=0}^d \theta_i^* E_i^*.
\end{equation}
Moreover
\begin{equation}   \label{productdual1}
A^*E_i^*\, =\, E_i^*A^*\, =\, \theta_i^* E_i^*  \ \ \ \ \ \ \ (0 \le i \le d).
\end{equation}
We now recall how $M$ and $M^*$ are related. 
By the definition of the distance matrices and dual idempotents, we have
\begin{equation}                                  \label{combcoeff}
E_h^* A_i E_j^* = 0\ \ \mbox{if and only if $p_{ij}^h=0$} \ \ \ \ \ \  
(0 \le h,i,j\le d).
\end{equation}
By \cite[Lemma 3.2]{Talg},
\begin{equation}                                  \label{propertiesdualbma}
E_h A_i^* E_j\,=\, 0\ \ \mbox{if and only if\ \ $q_{ij}^h=0$}  \ \ \ \ \ \ 
(0 \le h,i,j\le d).
\end{equation}

Let $T = T (x)$ denote the subalgebra of $\Mat_X\C$ generated by 
$M$ and $M^*$. 
We call  $T$ the {\de subconstituent algebra} or 
{\de Terwilliger algebra} of $\G$ with respect to $x$ \cite[p.~380]{Talg}.
By a {\de $T$-module} we mean a subspace $W \subseteq V$ such that 
$BW \subseteq W$ for all $B \in T$.
%
%
%
Let $W$ denote a $T$-module. Then $W$ is said to be {\de irreducible} 
whenever $W$ is nonzero and contains no $T$-modules other than 
$0$ and $W$.
%

%

We mention a special irreducible $T$-module. 
Let $\B{j}=\sum_{y \in X} {\hat y}$ denote the all 1's vector in $V$. 
Observe that $A_i {\hat x} = E^*_i \B{j}$ and 
$|X|E_i {\hat x} = A^*_i \B{j}$ for $0 \le i \le d$.
Therefore $M{\hat x} = M^* \B{j}$. 
Denote this common space by $V_0$ and note that $V_0$ is a $T$-module.
This $T$-module is irreducible by \cite[Lemma 3.6]{Talg}.
For $0 \le i \le d$ the vector $A_i{\hat x}$ is a basis for $E^*_iV_0$ 
and $E_i{\hat x}$ is a basis for $E_iV_0$. 
We call $V_0$ the {\it primary} $T$-module.

\section{The main result}

In this section we prove our characterization of $Q$-polynomial 
distance-regular graphs.

\begin{lemma}                \label{claim1} 
Let $\G$ denote a distance-regular graph with diameter $d \ge 3$. 
Let $E$ denote a nontrivial minimal idempotent of $\G$ that is TTR and 
let $\{ \theta_i^*\}_{i=0}^d$ be the corresponding dual eigenvalues. 
Let $\beta$ and $\gamma^*$ denote complex scalars such that
$\theta_{i-1}^*-\beta\theta_i^*+\theta_{i+1}^*=\gamma^*$ for 
$1 \leq i \leq d-1$.
Then the expression
\begin{equation}   \label{delta}
\theta^{*\,2}_{i-1}
-\beta\theta_{i-1}^* \theta_i^* +{\theta_i^*}^2
      -\gamma^*(\theta_{i-1}^*+\theta_i^*)
\end{equation} 
is independent of $i$ for $1\le i\le d$.
\end{lemma}

\proof 
Let $\delta_i^*$ denote the expression in (\ref{delta}).
For $1 \le i \le d-1$ the difference
$\delta_i^*-\delta_{i+1}^*$ is equal to
$$
(\theta_{i-1}^*-\theta_{i+1}^*)
(\theta_{i-1}^*-\beta\theta_i^*+\theta_{i+1}^*-\gamma^*),
$$
and is therefore 0. 
It follows that $\delta_i^*$ is
independent of $i$ for $1 \leq i \leq d$.  \qed

\begin{lemma}                \label{claim2}
Let $\G$ denote a distance-regular graph with diameter $d \ge 3$.
Let $E$ denote a nontrivial minimal idempotent of $\G$  that is TTR. 
Fix a vertex $x$ of $\G$ and let $A^*=A^*(x)$ denote the dual distance
matrix corresponding to $E$. Then
\begin{equation}                           \label{commutator}
0=[A^*,{A^*}^2 A -\beta A^*AA^*+A{A^*}^2
		-\gamma^*(AA^*+A^*A)-\delta^*A],
\end{equation}
where $\beta$ and $\gamma^*$ are from Lemma \ref{claim1} and 
$\delta^*$ it the common value of (\ref{delta}).
Here $[r,s]$ means $rs-sr$.
\end{lemma}

\proof 
Let $C$ denote the expression on the right in (\ref{commutator}).
We show that $C=0$. Observe
%
\begin{equation*}                
C = ICI = \biggl(\sum_{i=0}^d  E_i^*\biggr)\ C \ 
          \biggl(\sum_{j=0}^d  E_j^*\biggr)
  = \sum_{i=0}^d \sum_{j=0}^d \, E_i^* C E_j^*. 
\end{equation*}
To show $C=0$ it suffices to show $E_i^*CE_j^*=0$ 
for $0 \le i,j\leq d$.
For notational convenience define a polynomial $P$ in two
variables
$$
       P(u,v) = u^2 - \beta u v + v^2 - \gamma^*(u+v) - \delta^*.
$$
For $0 \leq i,j\leq d$ we have
$$
E^*_iCE^*_j = E^*_iAE^*_j\, 
P(\theta^*_i,\theta^*_j)\
(\theta^*_i-\theta^*_j)
$$
by (\ref{productdual1}), where $\{\theta^*_i\}_{i=0}^d$ are the dual 
eigenvalues corresponding to $E$. 
By (\ref{combcoeff}) we find  $E^*_iAE^*_j=0$ if $|i-j|>1$.
By Lemma \ref{claim1} we find $P(\theta^*_i,\theta^*_j) = 0$ if $|i-j|=1$.
Of course $\theta^*_i-\theta^*_j=0$ if $i=j$.
Therefore $E^*_iCE^*_j = 0$ as desired. We have now shown $C=0$.  \qed

\bigskip
In (\ref{propertiesdualbma}) we gave a characterization of a vanishing 
Krein parameter. We will need a variation on that result. 
We will obtain this variation after a few lemmas.
The first lemma follows directly from the definitions of 
$\hat{x}$ and the dual adjacency matrices.

\begin{lemma}                \label{claim}
Let $\G$ denote a distance-regular graph with diameter $d \ge 3$.
Let $E$ denote a nontrivial minimal idempotent of $\G$. 
Fix a vertex $x$ of $\G$ and let $A^*=A^*(x)$ denote the dual 
distance matrix corresponding to $E$. Then  
$
A^* \, v = |X| (E\,\hat{x}) \circ v 
$
for all $v\in V$. \hfill  \qed
\end{lemma}
%

\begin{lemma}                \label{claimne00}
Let $\G$ denote a distance-regular graph with diameter $d \ge 3$ and
let  $\{ E_i\}_{i=0}^d$ denote the minimal idempotents of $\G$.  
Fix $x \in X$.
%
%
Then for $0 \le h,i,j\le d$,
$$
E_h A^*_i E_j\, {\hat x} = q^h_{ij}\, E_h\, {\hat x},
$$
where $A^*_i = A^*_i(x)$.
\end{lemma}

\proof
From Lemma \ref{claim} we find
$
A^*_i E_j\, {\hat x} = |X|(E_i\, {\hat x}) \circ (E_j\, {\hat x}),
$
and observe that this equals $|X|(E_i \circ E_j)\, {\hat x}$. Now
$$
E_h A^*_i E_j\, {\hat x}
= |X| \, E_h (E_i \circ E_j)\, {\hat x}
= E_h   \sum_{\ell =0}^d q^\ell_{ij} E_\ell \, \hat x
= q^h_{ij} E_h \, {\hat x}. \ \ \ \ \qed
$$ 

 
\begin{corollary}                             \label{posledi}
Let $\G$ denote a distance-regular graph with diameter $d \ge 3$ and
let  $\{ E_i\}_{i=0}^d$ denote the minimal idempotents of $\G$.  
Fix $x \in X$ and let $V_0$ denote the primary module for $T(x)$.
Then the following {\rm (i), (ii)} are equivalent for $0 \le h,i,j\le d$.
\begin{description} \itemsep 35pt \parskip -30pt
\item[{\rm (i)}]  \hspace*{8pt}  $q_{ij}^h=0$.
\item[{\rm (ii)}] \hspace*{3pt}  $E_h A_i^* E_j$ 
                  vanishes on $V_0$, where $A^*_i = A^*_i(x)$.
\end{description}
\vspace*{-1mm}
Suppose {\rm (i), (ii)} fail. Then $E_h A^*_i E_jV_0 = E_h V_0$.
\end{corollary}  
\proof
\medskip \noindent
(i)$\implies$(ii) \ \ $E_h A_i^* E_j$ is zero
by (\ref{propertiesdualbma}) and hence vanishes on $V_0$.  

\noindent
(ii)$\implies$(i) \ \ 
Observe $\hat{x} \in V_0$ so $E_h A^*_i E_j\, {\hat x} = 0$.
Therefore $q^h_{ij} E_h\, {\hat x} = 0$ in view of Lemma \ref{claimne00}.
The vector $E_h\, {\hat x}$ is a basis for $E_hV_0$ so $E_h\, {\hat x}\ne 0$.
Thus $q^h_{ij}=0$. 
\medskip

\noindent
Suppose (i), (ii) fail. Then $E_h A^*_i E_jV_0$ is a nonzero subspace 
of the one-dimensional space $E_hV_0$ and is therefore equal to $E_hV_0$.
\qed

\bigskip
We are now ready to prove our main result.

\bigskip\noindent
{\it Proof of Theorem} \ref{main}.
First suppose that $E$ is $Q$-polynomial.
Then condition (i) holds by definition of a tail and
the definition of the $Q$-polynomial property,
condition (ii) holds by \cite[Theorem 8.1.2]{BCN}, cf. Leonard \cite{L},
and condition (iii) holds by \cite[Proposition 4.1.8]{BCN}.

\bigskip
To obtain the converse, assume that $E$ satisfies
(i)--(iii). We show $E$ is $Q$-polynomial. 
To do this we consider the representation diagram $\Delta_E$
from the introduction. We will show that $\Delta_E$ is a path.

Let $\{ E_i\}_{i=1}^d$ denote an ordering of the nontrivial minimal 
idempotents of $\G$ such that $E=E_1$.
Let $X$ denote the vertex set of $\G$. Fix $x \in X$ and let $A^*=A_1^*(x)$.
We first show that $\Delta_E$ is connected.
To do this we follow an argument given in 
\cite[Theorem 3.3]{Ter1}.
Suppose that $\Delta_E$ is not connected. Then there
exists a nonempty proper subset $S$ of $\{0,1,\dots,d\}$ 
such that $i$, $j$ are not adjacent in $\Delta_E$ for all 
$i\in S$ and $j \in \{0,1,\dots,d\}\backslash S$. 
Invoking (\ref{propertiesdualbma}) we find
$E_iA^*E_j = 0$ and $E_jA^*E_i = 0$ for $i \in S$ and
$j \in \{0,1,\ldots, d\} \backslash S$.
Define $F:={\sum_{i\in S} E_i}$ and observe
$$
A^* F
= IA^*F
= \biggl(\sum_{i=0}^d E_i\biggr) A^* F
= FA^* F.
$$
By a similar argument $F A^*= F A^* F$, so $A^*$ commutes with $F$.
Since $F \in M$ there exist complex scalars $\{\alpha_i\}_{i=0}^d$  such 
that $F = {\sum_{i=0}^d \alpha_i A_i}$. We have
\begin{equation}                                              \label{com}
0 = A^* F - F A^* = \sum_{i=1}^d \alpha_i (A^* A_i - A_i A^*).
\end{equation}
We claim that the matrices $\{ A^*A_i - A_iA^*\, |\, 1 \le i \le d \}$
are linearly independent.
To prove the claim, for $1 \leq i\leq d$  define $B_i = A^*A_i-A_iA^*$, 
and observe $B_i\, {\hat x} = (\theta^*_i-\theta^*_0)A_i\, {\hat x}$.
The vectors $\lbrace A_i\, {\hat x} \rbrace_{i=1}^d$ are linearly 
independent and
$\theta^*_i \not=\theta^*_0$ for $1 \leq i \leq d$ so the vectors
$ \lbrace B_i\, {\hat x} \rbrace_{i=1}^d$ are linearly independent.
Therefore the matrices $\lbrace B_i \rbrace_{i=1}^d$ are 
linearly independent and the claim is proved.
By the claim and (\ref{com}) we find $\alpha_i = 0$ for $1 \leq i \leq d$. 
%
%
%
Now $F = \alpha_0 I$. But $F^2=F$ so $\alpha^2_0 = \alpha_0$. 
Thus $\alpha_0=0$, in which case $S=\emptyset$, or $\alpha_0=1$, 
in which case $S=\{0,1,\ldots, d\}$.
In either case we have a contradiction so $\Delta_E$ is connected.

As we mentioned in the introduction, in $\Delta_E$
the vertex $0$ is adjacent to vertex $1$ and no other vertex of $\Delta_E$.
By the definition of a tail and since  $\Delta_E$ is connected, 
vertex $1$ is adjacent to vertex $0$ and exactly one other vertex in 
$\Delta_E$. To show that $\Delta_E$ is a path,
it suffices to show that each vertex in $\Delta_E$ is adjacent to
at most two other vertices in $\Delta_E$. We assume this is not
the case and obtain a contradiction. 
Let $v$ denote a vertex in $\Delta_E$ that is adjacent to more than 
two vertices of $\Delta_E$. Of all such vertices, we pick $v$ such that 
the distance to $0$ in $\Delta_E$ is minimal. Call this distance $i$.
Note that $2 \leq i \leq d-1$ by our above comments and the construction.
For notational convenience and without loss
of generality we may assume that the vertices of $\Delta_E$
are labelled such that for $1 \leq j \leq i-1$
the vertex $j$ is adjacent to $j-1$ and $j+1$ and no other vertex in $\Delta_E$.
By construction the chosen vertex $v$ is labelled $i$.
This vertex is adjacent to $i-1$ and at least two other vertices in $\Delta_E$.
Let $t$ denote a vertex in $\Delta_E$ other than $i-$1 
that is adjacent to vertex $i$. 
Since $E$ is TTR there exists $\beta \in \C$ such that
$\theta^*_{j-1}-\beta \theta^*_j + \theta^*_{j+1}$ is independent
of $j$ for $1 \leq j \leq d-1$. 
We claim that
\begin{equation}                                     \label{unique}
         \theta_t - (\beta +1)\theta_{i}  
                     + (\beta+1)\theta_{i-1}
                     -  \theta_{i-2}= 0.                   
\end{equation} 
%
To prove the claim we consider the equation (\ref{commutator}).
In that equation we expand the right-hand side to get
\begin{eqnarray*}
0=A^{*3}A -(\beta\! +\!1)(A^{*2}AA^*\!-\! A^*AA^{*2}) -AA^{*3}
        -\gamma^* (A^{*2}A\! -\! AA^{*2}) - \delta^* (A^*A\! -\! AA^*). 
\end{eqnarray*}
In this equation we multiply each term on the left by
$E_{i-2}$ and on the right by $E_t$. 
To help simplify the results we make some comments.
Using (\ref{prodAEi}) we find 
$E_{i-2}A^{*3}AE_t = E_{i-2}A^{*3}E_t \, \theta_t$.
Using (\ref{propertiesdualbma}) we find
$$
E_{i-2}A^{*3}E_t = E_{i-2} A^* \biggl(\sum_{r=0}^d E_r\biggr) A^*
  \biggl(\sum_{s=0}^d E_s\biggr) A^* E_t 
      = E_{i-2} A^* E_{i-1} A^* E_i A^* E_t.
$$
Similarly we calculate
\begin{eqnarray*}
E_{i-2} A^{*2} A A^*E_t &=& E_{i-2}A^*E_{i-1}A^*E_i A^*E_{t}\,\theta_i,\\
E_{i-2} A^*A A^{*2} E_t &=& E_{i-2}A^*E_{i-1}A^*E_i A^*E_{t}\,\theta_{i-1},\\
E_{i-2} A A^{*3}    E_t &=& E_{i-2}A^*E_{i-1}A^*E_i A^*E_{t}\,\theta_{i-2}
\end{eqnarray*}
and 
$$
\begin{array}{rccclrcl}
E_{i-2} A^{*2} A E_t&=&0, &  & \phantom{012345} &  E_{i-2} A A^{*2} E_t &=& 0, 
\\[2mm]
E_{i-2} A^* A E_t&=&0,    &  & \phantom{012345} &  E_{i-2} A A^* E_t &=& 0.    
\end{array}
$$
%
%
From these comments we find
\begin{equation}                                       \label{calculateh}      
0=E_{i-2}A^*E_{i-1}A^*E_i A^*E_t 
\bigl(\theta_t -(\beta+1)\theta_i +(\beta+1)\theta_{i-1} -\theta_{i-2}\bigr).
\end{equation}
We show that $E_{i-2}A^*E_{i-1}A^*E_iA^*E_t \ne 0$.
By the last assertion of Corollary \ref{posledi} and since
the sequence $(i-2,i-1,i,t)$ is a path in $\Delta_E$, we find
$E_iA^*E_tV_0=E_iV_0$ and 
$E_{i-1}A^*E_iV_0=E_{i-1}V_0$ and
$E_{i-2}A^*E_{i-1}V_0=E_{i-2}V_0$.
Therefore
$$
E_{i-2}A^*E_{i-1}A^*E_iA^*E_tV_0=E_{i-2}V_0.
$$
Observe $E_{i-2}V_0 \ne 0$ so 
$E_{i-2}A^*E_{i-1}A^*E_iA^*E_t \ne 0$, as desired.
By this and (\ref{calculateh}) we obtain (\ref{unique}). 
By (\ref{unique}) the scalar $\theta_t$ is uniquely determined.
The scalars $\theta_0,\ldots, \theta_d$ are mutually distinct
so $t$ is uniquely determined, for a contradiction. We have shown
that $\Delta_E$ is a path and therefore $E$ is $Q$-polynomial.
\qed

\section{Remarks}     \label{cx}

In this section we make some remarks concerning 
the three conditions in Theorem \ref{main}.
Throughout the section assume $\G$ is a distance-regular graph
with diameter $d \geq 3$ and 
eigenvalues $\theta_0 > \theta_1 > \cdots > \theta_d$.
Pick  a nontrivial minimal idempotent $E=E_j$ of $\G$ and let
$\lbrace \theta^*_i\rbrace_{i=0}^d$ denote the corresponding 
dual eigenvalue sequence. Abbreviate $\theta=\theta_j$.

\begin{note}                        \label{note:1}
\cite[pp.~142--143,\, 161]{BCN}
Pick an integer $i$ $(1 \leq i \leq d)$.
Then $\theta^*_i= \theta^*_0$ if and only if at least one of
the following holds.
\begin{description} \itemsep 35pt \parskip -30pt
\item[{\rm (i)}] \hspace*{3pt} $\G$ is bipartite, $i$ is even, and $j=d$,
\item[{\rm (ii)}] \hspace*{3pt} $\G$ is antipodal, $i=d$, and $j$ is even.
\end{description}
\end{note}

\begin{note}                         \label{note:2}
The graph $\G$ is imprimitive if and only if 
$\G$ is bipartite or antipodal \cite[Theorem 4.2.1]{BCN}. 
Thus if $\G$ is primitive then 
$\theta^*_i\not=\theta^*_0$ for $1 \leq i \leq d$.
\end{note}

\begin{lemma}                              \label{lem:3}
Assume that $\theta \not=-1$ and one of the following occurs:
\begin{description} \itemsep 35pt \parskip -30pt
\item[{\rm (i)}] \hspace*{3pt} $d=3$,
\item[{\rm (ii)}] \hspace*{3pt} $d=4$, $\G$ is antipodal, and $j$ is even,
\item[{\rm (iii)}] \hspace*{3pt} $d=5$, $\G$ is antipodal, and $j$ is even.
\end{description}
Then $E$ is TTR.
\end{lemma}

\proof
We have $(\theta^*_1-\theta^*_2)kb_1=(k-\theta)(1+\theta)$ 
by \cite[Lemma 2.2]{jkt},
so $\theta^*_1\not=\theta^*_2$. Define $\beta \in \C$ such that
$\beta+1 = (\theta^*_0-\theta^*_3)/(\theta^*_1-\theta^*_2)$.
By the construction 
$
\theta^*_{i-1}-\beta \theta^*_i+\theta^*_{i+1}
$
is independent of $i$ for $i=1,2$.
We are done in case (i), so assume we are in cases (ii) or (iii). 
By \cite[p.~142]{BCN} we have
$
\theta^*_i = \theta^*_{d-i}$ for $0 \leq i \leq d
$.
Therefore
$
\theta^*_{i-1}-\beta \theta^*_i+\theta^*_{i+1}
$ 
is independent of $i$ for $1 \leq i \leq d-1$. 
In other words $E$ is TTR.  \qed

\newpage
\begin{note}                               \label{note:4}
Referring to the conditions (i)--(iii) of Theorem \ref{main}, 
we show that no proper subset of (i)--(iii) implies that $E$ 
is $Q$-polynomial.
\begin{itemize}\itemsep -0pt
\leftskip -5mm
\item Conditions (ii), (iii) are not sufficient. 
      Assume $\G$ is the generalized hexagon of order $(2,1)$ 
      \cite[p.~200,\, 425]{BCN}.
      It is primitive with  diameter $d=3$ and 
      eigenvalues $4, 1+\sqrt{2}, 1-\sqrt{2}, -2$.
      Pick $j=1$. 
      Then $E$ satisfies condition (ii) by Lemma \ref{lem:3}
      and $E$ satisfies condition (iii) by Note \ref{note:2}.
      But $E$ does not satisfy (i) by \cite[p.~413,\, 425]{BCN}. 
      In particular $E$ is not $Q$-polynomial.

\item Conditions (i), (iii) are not sufficient.
      Assume $\G$ is the dodecahedron, which is antipodal with
      diameter $d=5$ and eigenvalues $3,\sqrt{5},1,0,-2,-\sqrt{5}$, 
		see \cite[p.~417]{BCN}.  Pick $j=1$. 
      Then $E$ satisfies condition (i) by \cite[p.~413,\, 417]{BCN} and 
      $E$ satisfies condition (iii) by  Note \ref{note:1}.
      By \cite[Lemma 2.2]{jkt} we find
      \ \ $
      \theta^*_0=3,
      \ \ \theta^*_1=\sqrt{5},
      \ \ \theta^*_2=1,
      \ \ \theta^*_3=-1,
      \ \ \theta^*_4=-\sqrt{5},
      \ \ \theta^*_5=-3
      $\ \ 
      and using this one verifies that $E$ does not satisfy condition (ii).
		In particular $E$ is not $Q$-polynomial.

\item Conditions (i), (ii) are not sufficient. 
      Assume $\G$ is the Wells graph \cite[p.~421]{BCN}, 
		which is antipodal with diameter $d=4$ and eigenvalues 
		$5,\sqrt{5},1,-\sqrt{5},-3$.  Pick $j=2$.  
      Then $E$ satisfies condition (i) by \cite[p.~413]{BCN} 
      and $E$ satisfies condition (ii) by Lemma \ref{lem:3}, 
      but $E$ does not satisfy condition (iii) by Note \ref{note:1}.
      In particular $E$ is not $Q$-polynomial.
\end{itemize}
\end{note}

\begin{note}                                    \label{note:5}
Assume $\G$ is bipartite, and consider the conditions
(i)--(iii) of Theorem \ref{main}.
If $E$ satisfies conditions (i), (iii) then $E$ is $Q$-polynomial
by \cite[proof of Theorem 5.4]{lang1}.
If $E$ satisfies conditions (ii), (iii) then $E$ is $Q$-polynomial 
by \cite[Theorem 10.5]{lang2}.
\end{note}

%
%
%



{\footnotesize
\begin{thebibliography}{99} 
\baselineskip 15pt
\itemsep -2pt

\bibitem{BI} 
Bannai, E. and T. Ito, {\it Algebraic Combinatorics I:
 Association Schemes}. Benjamin-Cummings Lecture Note Ser. 58,
The Benjamin/Cumming Publishing Company, Inc., London (1984).


\bibitem{BCN} 
A.\ E.\ Brouwer, A.\ M.\ Cohen and A.\ Neumaier, 
{\it Distance-Regular Graphs}, Springer-Verlag, Berlin, Heidelberg, 
1989.

\bibitem{delsarte} P. Delsarte,
An algebraic approach to the association schemes of coding theory,
{\em Philips Research Reports Suppl.} {\bf 10} (1973).

\bibitem{jkt}
A.\ Juri\v{s}i\'c, J.\ Koolen and P.\ Terwilliger,
Tight distance-regular graphs,
{\em J. Algebraic Combin.} {\bf 12} (2000), 163--197.

\bibitem{lang1} M. S. Lang,
Tails of bipartite distance-regular graphs,
{\em European J. Combin.} {\bf 23} (2002), 1015--1023.

\bibitem{lang2} M. S. Lang,
A new inequality for bipartite distance-regular graphs,
{\em J. Combin. Theory Ser. B} {\bf 90} (2004), 55--91.


\bibitem{L} D. A. Leonard, Orthogonal polynomials, duality and 
association schemes, {\em SIAM J. Math. Anal.} {\bf 13} (1982),
656-663.

\bibitem{pascasio3} A. A. Pascasio, 
A characterization of $Q$-polynomial distance-regular graphs, 
{\em Discrete Math.}  {\bf 308}  (2008), 3090--3096.

\bibitem{Talg} 
P.\ Terwilliger, 
The subconstituent algebra of an association scheme, I.,
{\em J. Algebraic Combin.} {\bf 1} (1992), 363--388.

\bibitem{Talg3} 
P.\ Terwilliger, 
The subconstituent algebra of an association scheme, III.,
{\em J. Algebraic Combin.} {\bf 2} (1993), 177--210.

 \bibitem{Ter1} 
P.\ Terwilliger, A new inequality for distance-regular graphs, 
{\em Discrete Math.} {\bf 137} (1995), 319--332.

\end{thebibliography}
}
\end{document}